%% file: arxiv.tex
\documentclass[onefignum,onetabnum]{siamsials251208}

\input{ex_shared}

\ifpdf
\hypersetup{
  pdftitle={Score-Based Generative Data Assimilation for Integrating Aggregated Surveillance Data into Agent-Based Models in Epidemic Tracking},
  pdfauthor={S. Liang, J. Hauck, M. Yang, A. Spannaus, H. Hanson, and G. Zhang}
}
\fi

\usepackage{xspace}
\usepackage{setspace}

\begin{document}

\maketitle

\begin{abstract}
Reliable epidemic monitoring often requires inferring regional infection burden and transmission heterogeneity from 
noisy, spatially aggregated, and potentially sparse surveillance data. Agent-based models (ABMs) are attractive for this task because they represent individual behavior, contact heterogeneity, and localized interventions, but these same features make them difficult to calibrate online. We develop a generative AI-based  data-assimilation (GenDA) framework for partially observed epidemic ABMs that estimates both the epidemic state and a heterogeneous parameter field while respecting the gap between observable macrostates and latent agent-level microstates. GenDA combines a training-free, score-based generative update for macrostate correction with a direct parameter update based on macrostate discrepancies, followed by a macro--micro reassignment step that restores consistency with the ABM.
In controlled and geographically explicit synthetic experiments, the framework recovers regional epidemic burden, dominant hotspot structures, and effective transmission heterogeneity from aggregated observations, while improving post-assimilation forecasts relative to state-only assimilation.
\end{abstract}

\begin{relevance}
Public-health surveillance rarely observes the full epidemic state. Instead, decision makers typically receive noisy, aggregated counts that obscure latent infection patterns, spatial heterogeneity, and changing transmission conditions. The proposed framework shows that, even under these challenging observational constraints, it is possible to extract reliable estimates of regional epidemic trends, dominant hotspot locations, and uncertainty summaries from heterogeneous ABMs without claiming exact recovery of individual trajectories. 
More broadly, the framework clarifies which epidemiologically meaningful quantities can be inferred reliably from aggregated surveillance: regional disease burden, dominant spatial risk patterns, and effective transmission and recovery heterogeneity, rather than unobservable individual trajectories. This distinction is important for interpreting and calibrating high-fidelity epidemic simulators used in outbreak monitoring and scenario planning.
\end{relevance}

\begin{mathcontent}
The methodology combines sequential Bayesian inference, score-based generative filtering, and joint state--parameter estimation for stochastic agent-based systems. Macroscopic epidemic counts are updated with an ensemble score filter that avoids linear-Gaussian assumptions, while an unknown parameter field is inferred through pseudo-observations defined in parameter space from macrostate discrepancies. An explicit macro--micro matching operator then maps updated aggregate states back to agent-level configurations. Numerical experiments quantify filtering accuracy, parameter recovery, and predictive behavior in both idealized and geographically explicit epidemic settings.
\end{mathcontent}
\vspace{-0.2cm}
\vspace{-0.2cm}

\input{0Introduction}
\input{1AgentBasedModel}
\input{3EnSF}

\input{4NumericalResults}
\input{5Conclusion}
\let\arxivstandardsection\section
\input{6Appendix}

\bibliographystyle{siamplain}
\bibliography{mybib}


\end{document}

%% file: ex_shared.tex
\usepackage{amsmath}
\usepackage{amssymb}
\usepackage{amsfonts}
\usepackage{graphicx}
\usepackage{subcaption}
\usepackage{color}
\usepackage{algorithmic}
\usepackage{bm}
\DeclareGraphicsExtensions{.pdf,.png,.jpg,.jpeg}

\usepackage{enumitem}
\setlist[enumerate]{leftmargin=.5in}
\setlist[itemize]{leftmargin=.5in}

\newsiamremark{remark}{Remark}
\newsiamremark{hypothesis}{Hypothesis}
\crefname{hypothesis}{Hypothesis}{Hypotheses}
\newsiamthm{claim}{Claim}
\newsiamremark{fact}{Fact}
\crefname{fact}{Fact}{Facts}

\makeatletter
\DeclareRobustCommand{\cev}[1]{%
  \mathpalette\do@cev{#1}%
}
\newcommand{\do@cev}[2]{%
  \fix@cev{#1}{+}%
  \reflectbox{$\m@th#1\vec{\reflectbox{$\fix@cev{#1}{-}\m@th#1#2\fix@cev{#1}{+}$}}$}%
  \fix@cev{#1}{-}%
}
\newcommand{\fix@cev}[2]{%
  \ifx#1\displaystyle
    \mkern#23mu
  \else
    \ifx#1\textstyle
      \mkern#23mu
    \else
      \ifx#1\scriptstyle
        \mkern#22mu
      \else
        \mkern#22mu
      \fi
    \fi
  \fi
}
\makeatother

\usepackage{amsopn}

\headers{GenDA for Inferring Regional Epidemic Dynamics}{S.~Liang, J.~Hauck, M.~Yang, A.~Spannaus, H.~Hanson, and G.~Zhang}

\title{Score-Based Generative Data Assimilation for Integrating Aggregated Surveillance Data into Agent-Based Models in Epidemic Tracking\thanks{This material was authored by UT-Battelle, LLC, under contract DE-AC05-00OR22725 with the US Department of Energy (DOE). The US government retains and the publisher, by accepting the article for publication, acknowledges that the US government retains a nonexclusive, paid-up, irrevocable, worldwide license to publish or reproduce the published form of this manuscript, or allow others to do so, for US government purposes. DOE will provide public access to these results of federally sponsored research in accordance with the DOE Public Access Plan.
\funding{This work was supported by funding from the Biopreparedness Research Virtual Environment (BRaVE) Program supported by
the U.S. Department of Energy, Office of Science, Advanced Scientific Computing Research and Biological and Environmental Research programs at Oak Ridge National Laboratory.}}}

\author{\hspace{1.4cm}Siming Liang\thanks{Computer Science and Mathematics Division, Oak Ridge National Laboratory, Oak Ridge, TN 37831, USA. }
\and Jacob Hauck\thanks{Mathematics and Statistics Department, Missouri University of Science and Technology, Rolla, MO 65401, USA.}
\and Minglei Yang\thanks{Fusion Energy Science Division, Oak Ridge National Laboratory, Oak Ridge, TN 37831, USA.}
\and Adam Spannaus\thanks{Computational Science and Engineering Division, Oak Ridge National Laboratory, Oak Ridge, TN 37831, USA.}
\and \vspace{0.05cm}\newline
Heidi Hanson\footnotemark[5]
\and Guannan Zhang\thanks{Corresponding author. Computer Science and Mathematics Division, Oak Ridge National Laboratory, Oak Ridge, TN 37831, USA. (\email{zhangg@ornl.gov})}}


%% file: 0Introduction.tex
\section{Introduction}
A recurring challenge in infectious-disease modeling is that the quantities most useful for real-time public-health decision making are often the hardest to observe directly. Surveillance systems typically provide noisy, delayed, and spatially aggregated counts of infections, hospitalizations, or test-confirmed cases, and reported case counts may also be affected by under-reporting and time-varying reporting rates \cite{spannaus2022inferring}. A full understanding of the latent epidemic burden, where high-risk regions are emerging, and how transmission conditions vary across populations and space is important for making informed decisions. Jointly estimating latent epidemic states and parameters from incomplete, aggregated surveillance data is therefore an ill-posed inverse problem: multiple agent-level configurations can be consistent with the same observations, so the underlying microstate is not uniquely identifiable. 

Agent-based models (ABMs) are attractive for this task because they represent individuals, their contacts, and their movement within their local environments. Relative to compartmental models, ABMs can encode household structure, heterogeneous mobility, behavioral adaptation, and geographically localized interventions, all of which matter when epidemic outcomes differ across neighborhoods, subpopulations, or activity hubs \cite{marshall2015formalizing,railsback_grimm_2019,zhang2025agent}. This flexibility has made ABMs increasingly important in outbreak analysis and scenario planning, from modular research platforms such as EpiPredict \cite{suer_epipredict_2024}, to national-scale systems such as UVA-EpiHiper \cite{bhattacharya_epihiper_2024}, to high-performance frameworks such as ENABLE \cite{spannaus_enable_2025}. The same features that make ABMs biologically and operationally appealing, however, also make them difficult to calibrate: the state is high-dimensional, the dynamics are stochastic and nonlinear, and many epidemiologically meaningful quantities are observed only after aggregation.

The resulting inferential problem has two coupled layers. At the micro level, the ABM evolves through agent states and interactions that are only partially or indirectly observable. At the macro level, public-health data report coarse summaries such as regional--level susceptible, infectious, and resistant counts or related surveillance indicators. The challenge is therefore not simply to estimate a hidden state, but to reconcile a latent agent-based system with the aggregated epidemic quantities that are actually measured. In realistic settings, exact reconstruction of the underlying microstate is neither possible nor necessary. What matters instead is whether the model can recover robust macro-level epidemic structure: regional burden, dominant hotspot patterns, and uncertainty ranges that remain informative for decision support \cite{gugole_uq_covid_abm_2021,knapp2025personalizing,swallow_emulation_abm_2022}.

Data assimilation provides a natural framework for this problem because it repeatedly merges model forecasts with incoming observations. Ensemble-based methods such as the ensemble Kalman filter and the local ensemble transform Kalman filter are computationally attractive and have been adapted to epidemiological ABMs \cite{clay2021real,hunt2007efficient,ward_enkf_abm_2016}. Their main limitation is that they rely on approximately linear-Gaussian structure, which is often violated in stochastic ABMs with nonlinear observation operators or multimodal forecast distributions \cite{bao2025nonlinear,xiong2025sensitivity}. Particle filters relax these assumptions but tend to suffer from degeneracy in high-dimensional systems \cite{gordon1993novel,sun2022analysis,tabataba2017epidemic,ternes2022data}. Direct state--parameter filters provide an alternative to augmenting parameters as state variables, but their adaptation to aggregated epidemic ABMs still requires a principled link between corrected macrostates and latent simulator parameters \cite{bao2023unified}. These difficulties become especially acute when the goal is not only to track epidemic state but also to infer heterogeneous transmission or recovery parameters online. Thus, the methodological gap addressed here is not the construction of another epidemic ABM, but the design of a sequential inference layer that can calibrate an existing stochastic ABM from sparse, aggregated surveillance data.

This paper develops GenDA, a generative-AI-based data-assimilation framework for partially observed epidemic ABMs. GenDA leverages the score-based filtering \cite{bao2024ensemble} for macrostate estimation and constructs a discrepancy-informed direct parameter update inspired by the Unified Filter \cite{bao2023unified}. Here ``generative AI'' refers to training-free score-based sampling, no neural network is trained in the filtering cycle. The central idea is to treat aggregate epidemic counts as the primary inferential target, use a non-Gaussian score-based filter to assimilate those quantities, and then transfer the resulting information into parameter space without requiring explicit parameter observations. A reassignment step restores consistency between the assimilated macrostate and the agent-level simulator. Framed this way, the method is not aimed at exact recovery of individual trajectories. It is aimed at extracting stable, epidemiologically meaningful information from complex ABMs under realistic surveillance constraints.

The main contributions are threefold. First, we formulate calibration of epidemic ABMs from aggregated surveillance as a coupled multiscale inference problem in which regional epidemic states, spatially heterogeneous parameters, and agent-level configurations must be updated consistently. 
Second, we develop GenDA, which introduces a discrepancy-informed parameter inference mechanism that transfers information from the macrostates correction to unobserved ABM parameters. Specifically, forecast–analysis discrepancies in the aggregated epidemic state are used to construct parameter-space pseudo-observations, enabling sequential estimation of spatially heterogeneous transmission and recovery parameters without requiring direct parameter observations. This parameter update is coupled with a non-Gaussian score-based macrostate filter and a macro–micro reassignment step that restores consistency between the assimilated aggregate state and the agent-level simulator.
Third, we validate the complete framework in controlled and geographically explicit ABMs and assess the contribution of parameter learning through direct comparison with the case where score filter is only used for macrostate update. Both methods improve state estimates during assimilation, but GenDA additionally reduces parameter error and maintains more accurate predictions after assimilation ends while recovering regional epidemic burden and hotspot structure from aggregated observations.

The rest of the paper is organized as follows. Section~\ref{sec:abm} defines the macroscopic and microscopic epidemic descriptions, the aggregated observation model, and the resulting model--data mismatch. Section~\ref{sec:method} presents the proposed GenDA framework and its macrostate estimation, parameter estimation, and macro--micro consistency modules. Section~\ref{sec:numericalexp} reports numerical experiments in idealized and geographically explicit settings, emphasizing epidemiological interpretation in addition to algorithmic performance. Section~\ref{sec:conclusion} concludes with implications, limitations, and future directions.

%% file: 1AgentBasedModel.tex
\section{Problem setting}\label{sec:abm}

This section introduces the modeling hierarchy that connects the quantities reported by public-health surveillance with the individual states evolved by an epidemic agent-based model (ABM). We first use a spatial susceptible--infectious--resistant (SIR) model as a macroscopic conceptual reference, then define the stochastic ABM used as the forecast model, and finally formulate the model--data mismatch that motivates the proposed data-assimilation framework.

\subsection{The PDE-based macroscopic epidemic model}

We begin with a spatially heterogeneous SIR model as a conceptual reference. This formulation provides an interpretable description of transmission, recovery, and spatial spread, but it is not used directly for simulation or inference in our framework. To incorporate spatial heterogeneity and mobility, a common extension of the classical ODE SIR model is a reaction--diffusion PDE formulation. Let $x \in \mathcal{D} \subset \mathbb{R}^d$ and $t \ge 0$, and define spatial densities
\[S = S(x,t),\quad I = I(x,t),\quad R = R(x,t),\quad N(x,t) = S(x,t) + I(x,t) + R(x,t).\]
A standard spatial SIR system is given by (see, e.g., \cite{keeling_rohani_epidemiology_2008,murray_mathematical_biology_2003})
\begin{equation}\label{eq:spdeSIR}
\begin{aligned}
\frac{\partial S}{\partial t}(x,t)
&= D_S\,\Delta S(x,t)
   - \beta(x,t)\,\frac{S(x,t)\,I(x,t)}{N(x,t)}
   + F_S(x,t)
\\[0.25em]
\frac{\partial I}{\partial t}(x,t)
&= D_I\,\Delta I(x,t)
   + \beta(x,t)\,\frac{S(x,t)\,I(x,t)}{N(x,t)}
   - \kappa(x,t)\,I(x,t)
   + F_I(x,t)
\\[0.25em]
\frac{\partial R}{\partial t}(x,t)
&= D_R\,\Delta R(x,t)
   + \kappa(x,t)\,I(x,t)
   + F_R(x,t),
\end{aligned}
\end{equation}
where $D_S,D_I,D_R\ge 0$ are diffusion coefficients that represent spatial mobility at an aggregate level. The rates $\beta(x,t)$ and $\kappa(x,t)$ allow spatially and temporally varying transmission and recovery, and $F_S,F_I,F_R$ collect external sources and sinks such as births, deaths, migration, vaccination, importation, or boundary inflow/outflow. Eq.~\eqref{eq:spdeSIR} summarizes population-level burden and spatial propagation, but its effective coefficients do not explicitly represent heterogeneous contacts, activity schedules, or discrete individual transitions.

\subsection{Agent-based microscopic epidemic model}

We therefore use an ABM as the dynamical forecast model. This choice does not assume that an ABM is universally more accurate than a compartmental model; rather, it permits mechanisms that are important in the applications considered here, including heterogeneous contacts, stochastic mobility, agent attributes, and geographically localized interventions \cite{railsback_grimm_2019,suer_epipredict_2024,zhang2025agent}. For notational simplicity, the ABM is defined directly in discrete time with $t=0,1,\ldots$ denoting the time step index. The complete agent-level microstate is defined by
\begin{equation}\label{eq:ABM_agents}
\mathbf A_t:=\bigl(\mathbf a_{t,1},\ldots,\mathbf a_{t,J}\bigr),
\end{equation}
where $\mathbf a_{t,j}$ denotes agent $j$ and contains its epidemiological state together with fixed or evolving attributes such as age group, location, compliance, and activity schedule. Let $\boldsymbol{\theta}_t:=(\boldsymbol \theta_{t,1},\ldots,\boldsymbol \theta_{t,J})$ denote the agent-indexed parameter field, where $\boldsymbol  \theta_{t,j}$ collects the epidemiological and behavioral parameters governing agent $j$, such as transmission and recovery parameters. This formulation permits fully agent-specific parameters, while an application may reduce the inferential dimension by sharing parameter values across any chosen spatial, demographic, behavioral, or other grouping. The stochastic ABM evolution is written generically as
\begin{equation}\label{eq:ABM_discrete}
\mathbf A_{t+1}=\mathcal F_t\bigl(\mathbf A_t,\boldsymbol{\theta}_t,\xi_t\bigr),
\end{equation}
where $\mathcal F_t$ contains the mobility, contact, behavioral, and epidemiological transition rules, and $\xi_t$ is a random variable that represents the stochasticity used by those rules. Thus, even for fixed $\mathbf A_t$ and $\boldsymbol{\theta}_t$, both $\mathbf A_{t+1}$ and any population statistic derived from it are random variables. Large-scale platforms demonstrate that this bottom-up representation can support heterogeneous, geographically detailed epidemic simulations, but they also highlight the computational and inferential difficulty of online calibration \cite{spannaus_enable_2025, bhattacharya_epihiper_2024}.

\subsection{The challenges in model-data integration in epidemic tracking}\label{sec:challenge}
Neither Eq.~\eqref{eq:spdeSIR} nor Eq.~\eqref{eq:ABM_discrete} is exact in practice. Epidemic forecasts are affected by uncertain and time-varying parameters, incomplete mobility and behavioral information, model-form error, and intrinsic stochasticity. These uncertainties can cause the simulated epidemic trajectory to drift from the evolving real system, which motivates calibration and the sequential integration of surveillance data \cite{gugole_uq_covid_abm_2021,kimpton_uq_abm_2024,li2017historymatching,spannaus2025robustuq,swallow_emulation_abm_2022}.

The central difficulty is that the ABM and the data generally live at different resolutions. 
To address this issue, we introduce a set of disjoint sub-regions of the spatial domain $\mathcal{D}$, i.e.,
\begin{equation}\label{eq:subregion}
    \mathcal{D} = \mathcal{D}_{1} \cup \cdots \cup  \mathcal{D}_K,
\end{equation}
over which surveillance statistics are reported. For the $k$-th sub-region, we define the epidemic macrostate
\begin{equation}\label{eq:macrostate}
\mathbf M_{t,k}:=\mathcal G_k(\mathbf A_t)
=\bigl(M_{S,t,k},M_{I,t,k},M_{R,t,k}\bigr)^\top,
\end{equation}
where the operator $\mathcal G_k$ aggregates all agents' information in the sub-region $\mathcal D_k$, and $M_{\cdot,t,k}$ is the sum of the epidemiological state in sub-region $k$, e.g., the $S, I, R$ states in Eq.~\eqref{eq:spdeSIR}. Because $\mathbf A_t$ is random, $\mathbf M_{t,k}$ is also a random vector. The corresponding surveillance model for the $k$-th sub-region is defined by
\begin{equation}\label{observations}
\mathbf Y_{t,k}=\mathcal H_k\bigl(\mathbf M_{t,k}\bigr)+\boldsymbol\varepsilon_{t,k},
\quad k=1,\ldots,K,
\end{equation}
where $\mathcal H_k$ may select observed compartments or represent an additional reporting transformation, and $\boldsymbol\varepsilon_{t,k}$ is observation and reporting error. We denote the collection of all sub-region level quantities as
\[
\mathbf M_t:=\{\mathbf M_{t,k}\}_{k=1}^{K},
\quad
\mathbf Y_t:=\{\mathbf Y_{t,k}\}_{k=1}^{K}.
\]

The inferential task is therefore to use the aggregated observations $\mathbf Y_{1:t}$ to update the random macrostates $\mathbf M_t$ and the parameter field $\boldsymbol{\theta}_t$, and then to modify the simulated agent states so that their sub-region  level counts agree with the updated macrostates. Prior work has shown both the value and the difficulty of assimilating coarse data into ABMs, particularly when parameters must be estimated and agent populations must be reconciled with aggregate corrections \cite{cocucci_enkf_epi_abm_2022,knapp2025personalizing,ternes2022data,ward_enkf_abm_2016}. The framework developed next addresses this model--data mismatch without treating the unobserved individual trajectory as uniquely identifiable.

%% file: 3EnSF.tex
\section{Score-based generative data assimilation for joint state and parameter estimation}\label{sec:method}
\begin{figure}[h!]
    \centering
\includegraphics[width=0.9\linewidth]{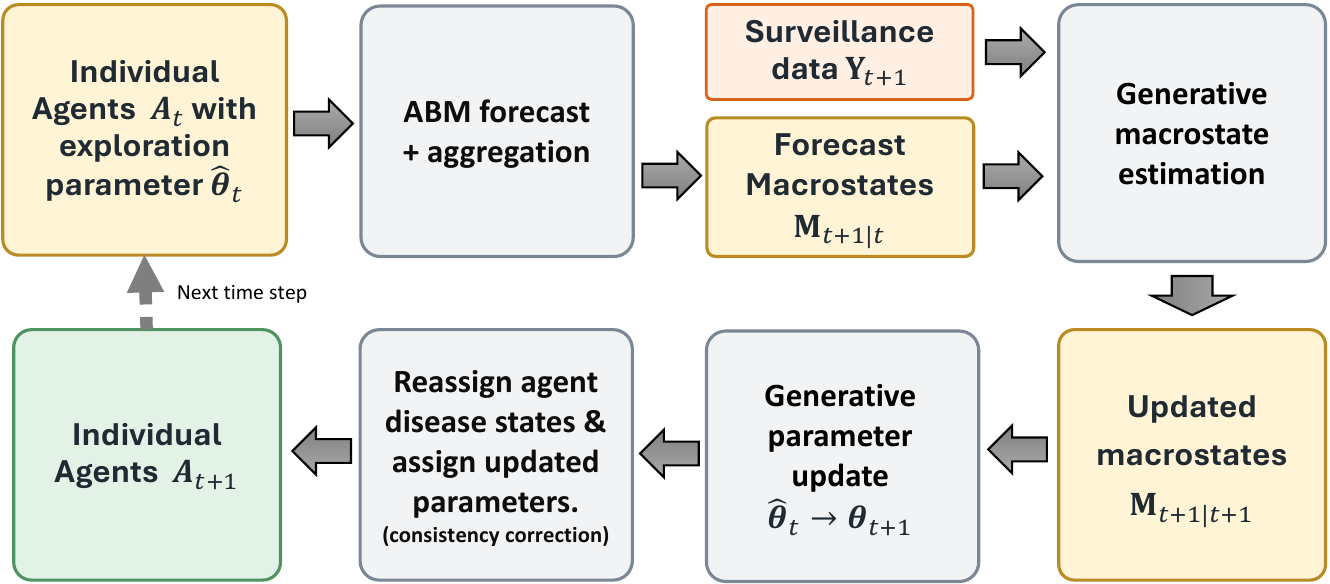}
    \caption{Overall GenDA forecast--update cycle. Agent state is propagated using exploratory parameters and aggregated over the observation sub-regions to obtain forecast macrostates. Surveillance observations are then used to update the macrostates. The agreement between each forecast and the updated macrostate determines how strongly its exploratory parameter values contribute to the parameter update. Agent disease states are subsequently reassigned locally to match the member-specific updated macrostates, while the updated parameters $\boldsymbol{\theta}_{t+1}$ are assigned to the corresponding agents or parameter-sharing groups. The reconstructed agents and updated parameters initialize the next forecast cycle.}
\label{fig:genda_overall}
\end{figure}

This section introduces the proposed GenDA framework designed to jointly estimate epidemic macrostates and unknown model parameters in stochastic ABMs using only aggregated observations. As summarized in Figure~\ref{fig:genda_overall}, each assimilation cycle has three components: (I) nonlinear filtering of the macrostates on the observation sub-regions, (II) sequential inference of the ABM parameter field, and (III) enforcement of macro--micro consistency by modifying agent disease states. Sections~\ref{sec:score_filter}--\ref{sec:Micro--macro} describe these components, and Section~\ref{sec:DAfilter} summarizes their coupling.

\subsection{Score-based nonlinear filtering for macrostate estimation} \label{sec:score_filter} 
We begin with the macrostate update. From simulation time $t$ to $t+1$, the agent state, its aggregated realization, and the observation model satisfy
\begin{equation}\label{eq:state_obs_system}
\begin{aligned}
\text{State:}\qquad  &\mathbf A_{t+1} 
=\mathcal F_t (\mathbf A_{t} ,\widehat{\boldsymbol \theta}_{t} ,\xi_t ),\\
&\mathbf M_{t+1,k} 
=\mathcal G_k\!\left(\mathbf A_{t+1} \right),\\
\text{Observation:}\qquad  &\mathbf Y_{t+1,k}
=\mathcal H_k\!\left(\mathbf M_{t+1,k}\right)+\boldsymbol\varepsilon_{t+1,k},\quad k=1,\ldots,K,
\end{aligned}
\end{equation}
where $\widehat{\boldsymbol \theta}_{t} $ is the exploration parameter field used for the forecast (details will be provided in Section~\ref{sec:directfilter}). The objective is to characterize the filtering distribution $p({\mathbf M_{t+1}|\mathbf Y_{1:t+1}})$ of $\mathbf M_{t+1}$ conditioned on the observation history $\mathbf Y_{1:t+1}:=\{\mathbf Y_1,\ldots,\mathbf Y_{t+1}\}$. Sequential Bayesian filtering consists of a prediction step followed by an update step. In the prediction step, propagation through the ABM and aggregation operators gives
\begin{equation}\label{eq:filtering-prior-ck}
p({\mathbf M_{t+1}|\mathbf Y_{1:t}})
=\int p({\mathbf M_{t+1}|\mathbf M_t})
p({\mathbf M_t|\mathbf Y_{1:t}})\,d\mathbf M_t.
\end{equation}
After the observation $\mathbf Y_{t+1}$ becomes available, Bayes' theorem gives
\begin{equation}\label{eq:filtering-bayesian-update}
\underbrace{p({\mathbf M_{t+1} \mid \mathbf Y_{1:t+1}})}_{\rm posterior}
\propto
\underbrace{p({\mathbf M_{t+1}|\mathbf Y_{1:t}})}_{\rm prior}
\underbrace{p({\mathbf Y_{t+1}|\mathbf M_{t+1}})}_{\rm likelihood}.
\end{equation}
For Gaussian observation error $\boldsymbol\varepsilon_{t+1}\sim\mathcal N(0,\boldsymbol\Sigma_{\varepsilon})$ and the stacked observation operator $\mathcal H$, the likelihood is
\begin{equation}\label{eq:likelihood}
p({\mathbf Y_{t+1}\mid\mathbf M_{t+1}})
\propto
\exp\left[-\frac12
\bigl(\mathbf Y_{t+1}-\mathcal H(\mathbf M_{t+1})\bigr)^\top
\boldsymbol\Sigma_{\varepsilon}^{-1}
\bigl(\mathbf Y_{t+1}-\mathcal H(\mathbf M_{t+1})\bigr)\right].
\end{equation}

The score-based nonlinear filtering method \cite{bao2024ensemble,bao2024score} provides a training-free generative mechanism to perform the update step Eq.~\eqref{eq:filtering-bayesian-update} without relying on linear-Gaussian assumptions. It represents the prior distribution using an ensemble-based score approximation and constructs samples from the posterior by solving a pair of forward and reverse SDEs over a synthetic diffusion time $\tau\in[0,1]$:
\begin{equation}\label{eq:sdes}
\begin{aligned}
   & \text{Forward:}\quad {\rm d} \mathbf{Z}_{t+1,\tau} = b_\tau \mathbf{Z}_{t+1,\tau}\, \rm d\tau + \sigma_\tau\, \rm d \mathbf{W}_\tau,\\[2pt]
   & \text{Reverse:}\quad {\rm d} \mathbf{Z}_{t+1,\tau} = \big[b_\tau \mathbf{Z}_{t+1,\tau} - \sigma_\tau^2 \mathbf{S}_{t+1|t}(\mathbf{Z}_{t+1,\tau}, \tau)\big] \rm d\tau + \sigma_\tau\, \rm d\cev{\mathbf{W}}_\tau,
\end{aligned}
\end{equation}
where $\mathbf{Z}_{t+1,\tau}$ denotes the SDE state at filtering step $t+1$ and diffusion time $\tau$, $\mathbf{S}_{t+1\mid t}(\mathbf{Z}_{t+1,\tau},\tau)$ is the score function of the prior distribution in Eq.~\eqref{eq:filtering-bayesian-update}, and $\mathbf{W}_\tau$, $\cev{\mathbf{W}}_\tau$ are the forward and backward Brownian motions. Even though there are multiple choices for the drift and diffusion coefficients in Eq.~\eqref{eq:sdes}, we use the following definition in this work, i.e.,
\begin{equation}\label{eq:cof}
b_\tau = \frac{\rm d\log \alpha_\tau}{\rm d\tau}, \qquad 
\sigma_\tau^2 = \frac{\rm d\omega_\tau^2}{\rm d\tau} - 2\,\frac{\rm d\log \alpha_\tau}{\rm d\tau}\,\omega_\tau^2,
\end{equation}
with $\alpha_\tau = 1 - \tau$ and $\omega_\tau^2 = \tau$. By initializing the forward SDE with samples from the prior filtering distribution at $\tau=0$, the forward transport maps that distribution to the standard Gaussian distribution $\mathcal N(0,\mathbf I_d)$ at $\tau=1$. Conversely, the reverse SDE maps standard Gaussian samples back to samples from the filtering distribution using the score $\mathbf{S}_{t+1\mid t}(z_{t+1,\tau},\tau)$.

\begin{figure}[!tb]
    \centering
    \includegraphics[width=0.99 \linewidth]{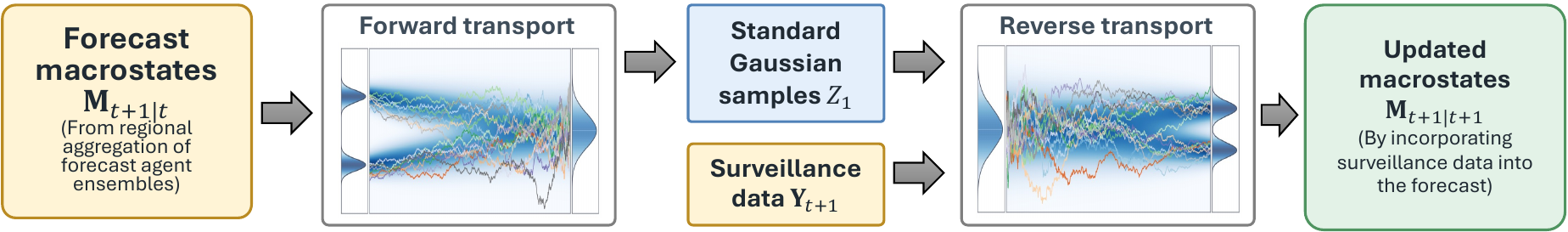}
    \caption{Illustration of the score-based filtering for macrostate update. The forecast macrostate ensemble obtained from the agents in Figure~\ref{fig:genda_overall} is transported via the forward diffusion model to a standard Gaussian reference distribution. During the reverse diffusion transport, the prior score is combined with likelihood information from the surveillance observations to generate the updated macrostate ensemble. The colored curves represent generative transport paths in diffusion time and do not represent epidemic evolution. The Gaussian distribution is a temporary computational reference used to facilitate posterior sampling.}
    \label{fig:flow_macrostate}
\end{figure}

Based on the diffusion model framework, we define the prior and posterior filtering distributions through $L$ ensembles as the following:
\begin{equation}\label{eq:prior_sample}
\mathcal P_{t+1}^{\rm prior}
:=\left\{\mathbf M_{t+1\mid t,l} \right\}_{l=1}^{L}, \quad \mathcal P_{t+1}^{\rm posterior}
:=\left\{\mathbf M_{t+1\mid t+1,l} \right\}_{l=1}^{L}. 
\end{equation}
As the drift term of the forward SDE in Eq.~\eqref{eq:sdes} is linear, a closed-form expression for the score function can be derived and approximated via a training-free Monte Carlo estimator:
\begin{equation}
\mathbf{S}_{t+1|t}(\mathbf{Z}, \tau) \;\approx\; 
\hat{\mathbf{S}}_{t+1|t}(\mathbf{Z}, \tau) 
:= \sum_{i=1}^{N_{\rm batch}} 
    -\,\frac{\mathbf{Z} - \alpha_\tau \mathbf M_{t+1|t,i}}{\omega_\tau^2}
\,\hat{\zeta}_\tau\left(\mathbf{Z}, \mathbf M_{t+1|t,i}\right),
\label{eq:score_mc}
\end{equation}
where $\{\mathbf M_{t+1\mid t,i}\}_{i=1}^{N_{\rm batch}}$ is a minibatch of samples from the forecast ensemble, and $\hat{\zeta}_\tau$ is a Monte Carlo approximation of the transition weights,
\begin{equation}
\hat{\zeta}_\tau\left(\mathbf{Z}, \mathbf M_{t+1|t,i}\right)
:= 
\frac{q\left(\mathbf{Z} \mid \mathbf M_{t+1|t,i}\right)}{
\sum_{i'=1}^{N_{\rm batch}} q\left(\mathbf{Z} \mid \mathbf M_{t+1\mid t,{i'}}\right)},
\label{eq:weights_mc}
\end{equation}
where $q(\mathbf{Z}| \cdot)$ is the transition density associated with the forward SDE.
To incorporate observations, the score-based nonlinear filtering combines the prior score with the likelihood gradient to form a posterior score. The Monte Carlo estimate of the posterior score is then defined as
\begin{equation}
\mathbf{S}_{t+1|t+1}(\mathbf{Z}, \tau) 
:= \mathbf{S}_{t+1|t}(\mathbf{Z}, \tau) 
+ h(\tau)\, \nabla_{\mathbf{Z}} \log p(\mathbf {Y}_{t+1} \mid \mathbf{Z}),
\quad \tau \in [0, 1],
\label{eq:score_update}
\end{equation}
where $h(\tau)$ is a monotone damping function with $h(0)=1$ and $h(1)=0$, which gradually injects likelihood information during the reverse diffusion in accordance with Bayes' rule. Finally, the updated ensemble $\left\{\mathbf M_{t+1\mid t+1,l} \right\}_{l=1}^{L}$ is obtained by sampling $\mathbf{Z}_{t+1,1} \sim\mathcal{N}(0,\mathbf I_d)$ and integrating the reverse SDE in Eq.~\eqref{eq:sdes} with the approximate posterior score in Eq.~\eqref{eq:score_update}.
Figure~\ref{fig:flow_macrostate} summarizes the forward transport, incorporation of surveillance information, and reverse transport that produce the updated macrostate ensemble.


\subsection{Parameter estimation through the data-informed macrostate} \label{sec:directfilter}
Epidemiological parameters such as transmission and recovery rates are generally uncertain, and their effective values may vary across individuals, populations, locations, or time \cite{cocucci_enkf_epi_abm_2022,kim2021automatic,swallow_emulation_abm_2022}. These parameters are not directly observed by epidemiological surveillance, so a standard parameter-observation operator is unavailable. Following the direct-filter principle in \cite{bao2023unified}, we instead construct parameter pseudo-observations from discrepancies between paired ABM forecasts and the data-informed macrostate update. The key insight is that parameters influence the epidemic dynamics only through the model evolution. Ensemble members whose parameter values are closer to the truth tend to produce state predictions (prior) that are more consistent with data-informed state estimates (posterior mean). By exploiting this relationship, we transform the updated macrostate obtained through Section~\ref{sec:score_filter} into an implicit observation of the parameters. Sequential Bayesian inference can then be performed in parameter space using the score-based filter.

\remark The proposed parameter estimation method can handle any application-specific parameter components or sharing groups such as spatial, demographic, behavioral, or other sharing structures without changing the filtering procedure. For notational simplicity, we suppress a separate parameter-group index and present the case in which the parameter partition coincides with the observation partition. Thus, each sub-region $\mathcal D_k$ is associated with one parameter vector shared within that sub-region.

\subsubsection{Constructing pseudo-observations in parameter space}
To calibrate the unknown parameters, we first introduce a pseudo-evolution model for adaptive parameter exploration. At the start of a forecast cycle, the parameter ensemble is given an exploration spread which can be interpolated as the forward model of the parameter field. Conceptually,
\begin{equation}\label{parameters}
\widehat{\boldsymbol\theta}_{t} 
=\boldsymbol\theta_{t} +\boldsymbol\eta_{t} ,
\qquad
\boldsymbol\eta_{t} \sim
\mathcal N\!\left(0,\gamma_{t}^2\mathbf I_{d_\theta}\right),
\end{equation}
where the scalar $\gamma_{t}>0$ is the prescribed exploration standard deviation and $d_\theta$ is the dimension of the stacked parameter vector. The schedule may decrease over time to move from coarse exploration to refined estimation.

Following the ensemble notation in Eq.~\eqref{eq:prior_sample}, we extend the notation  to include the parameter in each ensemble. The prior ensemble is defined by
\begin{equation}\label{eq:paired_parameter_forecast}
\left\{\left(\mathbf M_{t+1\mid t,l} ,
\widehat{\boldsymbol\theta}_{t,l} \right),\;\; l = 1, \ldots, L \right\}.
\end{equation}
To avoid outliers while retaining precise control of the exploration range in the numerical implementation, rather than adding independent perturbations to the parameter ensemble mean, we rescale the componentwise spread of the existing ensemble to match $\gamma_{t}$:
\begin{equation}\label{eq:propagateparamenter}
\widehat{\boldsymbol\theta}_{t,l} 
=\left(\boldsymbol\theta_{t,l} -\overline{\boldsymbol\theta}_{t}\right)
\cdot\frac{\gamma_{t}}
{\operatorname{std}\!\left(\{\boldsymbol\theta_{t,l} \}_{l=1}^{L}\right)}
+\overline{\boldsymbol\theta}_{t},
\qquad
\overline{\boldsymbol\theta}_{t}=\frac1{L}\sum_{l=1}^{L}\boldsymbol\theta_{t,l} ,
\end{equation}

Then, we construct the pseudo parameter observations. The likelihood weight is
\begin{equation}\label{eq:parameter_likelihood}
W_{t+1,l}  \propto
\exp\!\left[-\frac12
\bigl(\mathbf d_{t+1,l } \bigr)^\top
(\nu^2_{t+1 ,l} \mathbf{I}_{d_\theta})^{-1}
\mathbf d_{t+1,l } \right], \qquad \mathbf  d_{t+1,l} = \mathbf M_{t+1\mid t,l} -\overline{\mathbf M}_{t+1|t+1},
\end{equation}
where $\overline{\mathbf M}_{t+1|t+1}$ is the updated ensemble mean of the macrostates defined by 
\begin{equation}\label{eq:analysis_macrostate_mean}
\overline{\mathbf M}_{t+1\mid t+1}
:=\frac{1}{L}\sum_{l=1}^{L}\mathbf M_{t+1\mid t+1,l},
\end{equation}
%
%
and $\mathbf  d_{t+1,l}$ calculates the macrostate forecast--update discrepancy corresponding to each parameter. The scaling factor $\nu_{t+1,l}$ in Eq.~\eqref{eq:parameter_likelihood} is defined by 
\begin{equation}\label{eq:parameter_omega}
\nu_{t+1,l}
=\frac1{L}\sum_{l=1}^{L}
\left\|\mathbf d_{t+1,l} \right\|_1,
\end{equation}
which can also be a user-chosen or empirically estimated standard deviation that controls how discrepancies are weighted. Details on numerical implementation for $\nu_{t+1,l}$ and additional choices are provided in the supplementary material. The pseudo-observation in parameter space is then defined as the likelihood-weighted average conditioned on $\overline{\mathbf M}_{t+1|t+1}$, using a normalized weighted sum over the parameter ensemble:
\begin{equation}\label{eq:parameter_obs}
\begin{aligned}
\boldsymbol\theta_{t+1}^{\mathrm{obs}}
&:=\sum_{l=1}^{L}\widehat{\boldsymbol\theta}_{t,l} 
\frac{W_{t+1,l} }{\sum_{l'=1}^{L}W_{t+1,l'}}.
\end{aligned}
\end{equation}
where the error in the parameter observation is assumed to follow $\mathcal{N} (0,\gamma_{t+1}^2)$.

\subsubsection{Sequential Bayesian inference for the parameter field} Now we utilize the pseudo-observations $\boldsymbol\theta_{t+1}^{\mathrm{obs}}$ obtained above to perform sequential Bayesian inference of the parameter field. 
The prior parameter distribution is represented by $\{\widehat{\boldsymbol\theta}_{t,l} \}_{l=1}^{L}$. Bayes' theorem gives
\begin{equation}\label{Bayes:para_revised}
\underbrace{p\!\left(\widehat{\boldsymbol\theta}_{t}\mid
\overline{\mathbf M}_{1:t+1}\right)}_{\rm posterior}
\propto
\underbrace{p\!\left(\widehat{\boldsymbol\theta}_{t}\mid
\overline{\mathbf M}_{1:t}\right)}_{\rm prior}
\underbrace{p\!\left(\boldsymbol\theta_{t+1}^{\mathrm{obs}}\mid
\widehat{\boldsymbol\theta}_{t}\right)}_{\rm pseudo\text{-}likelihood}.
\end{equation}
We model the pseudo-observation error with covariance $\gamma_{t+1}^2\mathbf I_{d_\theta}$:
\begin{equation}\label{Bayes:paraobslike_revised}
p\!\left(\boldsymbol\theta_{t+1}^{\mathrm{obs}}\mid
\widehat{\boldsymbol\theta}_{t}\right)
\propto
\exp\!\left[-\frac12
\bigl(\widehat{\boldsymbol\theta}_{t}-\boldsymbol\theta_{t+1}^{\mathrm{obs}}\bigr)^\top
\left(\gamma_{t+1}^2\mathbf I_{d_\theta}\right)^{-1}
\bigl(\widehat{\boldsymbol\theta}_{t}-\boldsymbol\theta_{t+1}^{\mathrm{obs}}\bigr)\right].
\end{equation}

\begin{figure}[h!]
    \centering
\includegraphics[width=0.95\linewidth]{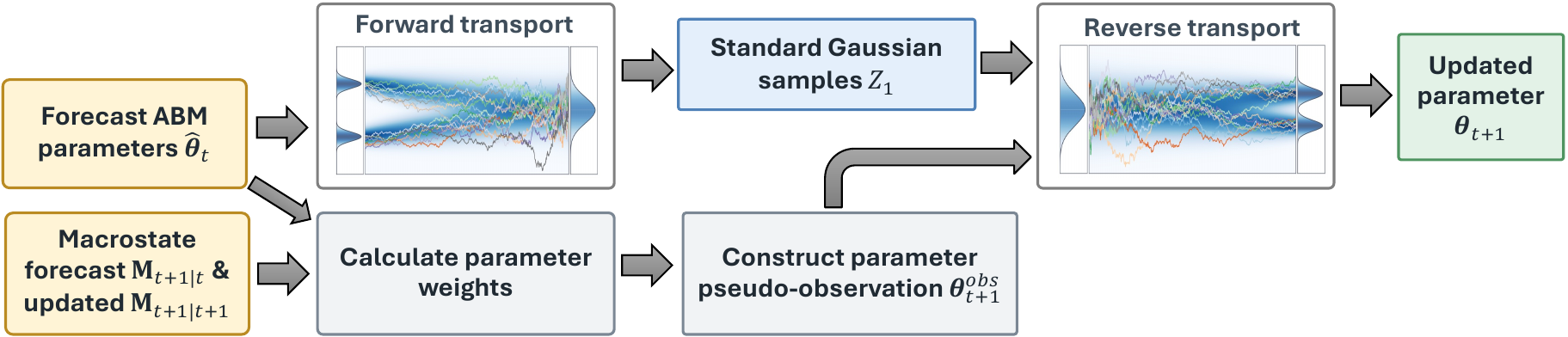}
    \caption{Discrepancy-informed parameter update. Each exploratory parameter realization is paired with the forecast macrostate that it produces. The forecast is compared with the updated macrostate mean obtained from Figure~\ref{fig:flow_macrostate}; ensemble members whose forecasts agree more closely with the updated macrostate receive larger likelihood weights. The weighted average of the paired exploratory parameters defines a parameter pseudo-observation, which is assimilated by the parameter-space score-based filter to obtain the updated parameter.}
\label{fig:flow_parameter}
\end{figure}

To perform the update, we combine the likelihood gradient from Eq.~\eqref{Bayes:paraobslike_revised} with the prior parameter score, approximated using the training-free Monte Carlo estimator in Eq.~\eqref{eq:score_mc}. Solving the reverse diffusion SDE then produces the updated parameter ensemble
\begin{equation}\label{eq:parameter_posterior_ensemble}
\left\{\boldsymbol\theta_{t+1,l} \right\}_{l=1}^{L}
\sim p\!\left(\widehat{\boldsymbol\theta}_{t}\mid
\overline{\mathbf M}_{1:t+1}\right),
\qquad
\overline{\boldsymbol\theta}_{t+1}
=\frac1{L}\sum_{l=1}^{L}\boldsymbol\theta_{t+1,l} .
\end{equation}
Figure~\ref{fig:flow_parameter} summarizes the construction of discrepancy-based likelihood weights, the parameter pseudo-observation, and the score-based parameter update.

\subsection{Macro--micro consistency enforcement}\label{sec:Micro--macro}
After each data assimilation step, the epidemic macrostate has been updated using aggregated observations. However, the agent-based simulator evolves at the level of individual agents, and the updated macrostate generally does not correspond to the current configuration of agent-level disease states. To close the data assimilation--ABM cycle, it is therefore necessary to reconcile the microstate of agents with the assimilated macrostate on the observation sub-regions.

\begin{figure}[h!]
    \centering
    \includegraphics[width=0.96\linewidth]{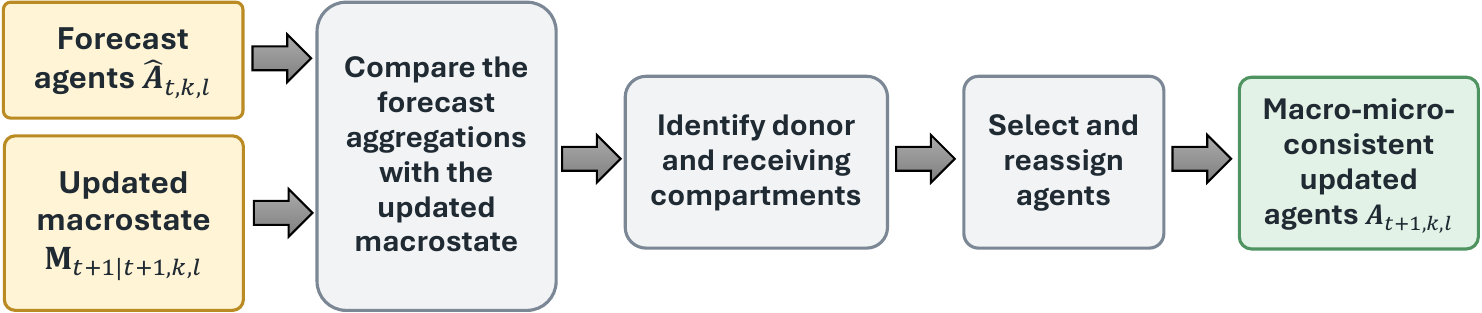}
    \caption{Macro--micro reassignment on one observation sub-region $k$ of one ensemble $l$. Given updated macrostate from Figure~\ref{fig:flow_macrostate} and forecast agents from Figure~\ref{fig:genda_overall}, identify the mismatch and reassign agents to ensure macro-micro consistency.}
    \label{fig:flow_reassignment}
\end{figure}

We focus on one ensemble member within one observation sub-region $\mathcal D_k$. First, we compare the forecast agent states with the updated macrostate to identify the updated-minus-forecast discrepancy. Then, agents are randomly selected without replacement from compartments with a surplus and reassigned to compartments with a deficit. Only their disease states are changed, their locations and all other attributes remain unchanged. This operation is applied locally to every sub-region and to every analysis ensemble member, using the member-specific updated macrostate rather than the posterior mean. If several compartments have surpluses and deficits, transfers continue without replacement until all target counts are met. The implementation used in the experiments samples donors uniformly. When agent-level risk information is available, a probability-weighted alternative can prioritize agents using infection risk or contact-history information. Formal donor--recipient notation and the weighted alternative are provided in Supplementary. The required output condition is
\begin{equation}\label{eq:macro_micro_consistency}
\mathcal G_k\!\left(\mathbf A_{t+1,l} \right)
=\mathbf M_{t+1\mid t+1,k,l} ,
\qquad l=1,\ldots,L,\quad k=1,\ldots,K.
\end{equation}
In both cases, the resulting agent configuration satisfies exact micro--macro consistency within each spatial cell after every assimilation step. This reassignment step completes the data assimilation--ABM cycle by ensuring that subsequent forward simulations begin from an agent-level state that is fully consistent with the assimilated macroscopic information.
Figure~\ref{fig:flow_reassignment} shows the module interface and the before/after example. Feasibility projection of the updated macrostate to nonnegative integer SIR counts is described in the supplementary material.

\subsection{The complete GenDA workflow}\label{sec:DAfilter}

Algorithm~\ref{alg:genda} combines the three modules into one sequential forecast--update cycle. The two outputs required for the next ABM forecast are the macro--micro-consistent agent ensemble from Section~\ref{sec:Micro--macro} and the updated parameter ensemble from Section~\ref{sec:directfilter}.

\begin{algorithm}
\caption{GenDA framework for joint macrostate and parameter estimation}
\label{alg:genda}
\setstretch{1.2}
\begin{algorithmic}[1]
\STATE{\textbf{Input}: ABM maps $\{\mathcal F_t\}$, aggregation maps $\{\mathcal G_k\}_{k=1}^{K}$, observation maps $\{\mathcal H_k\}_{k=1}^{K}$,  observations $\{\mathbf Y_t\}$, and exploration schedule $\{\gamma_t\}$}
\STATE{Initialize $\{\mathbf A_{0,l} \}_{l=1}^{L}$ and $\{\boldsymbol\theta_{0,l} \}_{l=1}^{L}$}
\FOR{$t=0,\ldots,T-1$}
\STATE{Construct exploration parameters $\{\widehat{\boldsymbol\theta}_{t,l} \}_{l=1}^{L}$ using Eq.~\eqref{eq:propagateparamenter}}
\STATE{Propagate the ABMs to obtain $\{\widehat{\mathbf A}_{t+1,l} \}_{l=1}^{L}$ using Eq.~\eqref{eq:ABM_discrete}}
\STATE{Aggregate the agents to get $\{\mathbf M_{t+1\mid t,l} \}_{l=1}^{L}$ and form the paired ensemble in Eq.~\eqref{eq:paired_parameter_forecast}}
\STATE{Assimilate $\mathbf Y_{t+1}$ with the score-based filter to obtain $\{\mathbf M_{t+1\mid t+1,l} \}_{l=1}^{L}$ and $\overline{\mathbf M}_{t+1\mid t+1}$}
\STATE{Construct $\boldsymbol\theta_{t+1}^{\mathrm{obs}}$ using Eqs.~\eqref{eq:parameter_likelihood}--\eqref{eq:parameter_obs}}
\STATE{Apply the parameter-space score-based filter update to obtain $\{\boldsymbol\theta_{t+1,l} \}_{l=1}^{L}$}
\STATE{Reassign agent disease states to obtain $\{\mathbf A_{t+1,l} \}_{l=1}^{L}$ satisfying Eq.~\eqref{eq:macro_micro_consistency}}
\ENDFOR
\STATE{\textbf{Output}: sequential ensembles of updated macrostates, parameters, and macro--micro-consistent ABM realizations}
\end{algorithmic}
\end{algorithm}

%% file: 4NumericalResults.tex
\section{Numerical experiments: recovering epidemic burden, hotspot structure, and transmission heterogeneity} \label{sec:numericalexp}

In this section we present numerical experiments to demonstrate the performance of the proposed framework for (i) accurate {state} estimation and (ii) robust {online parameter} estimation in stochastic, agent-based epidemic models. All experiments are implemented using the \texttt{Mesa} framework \cite{ter2025mesa}, a flexible Python toolkit for rapid ABM development. Although a pure-Python implementation is not intended for extreme-scale simulations, it provides a transparent and modular interface that is well suited for validating data assimilation ideas and testing algorithmic components. For geographically explicit experiments, we also use the extension \texttt{Mesa-Geo} \cite{wang2022mesa-geo}, which supports agents evolving on real map geometries. In both numerical examples, the permitted epidemiological transitions are $S\to I$, $I\to S$, and $I\to R$. Thus, an infected agent who clears the infection may return to the susceptible state and later be reinfected, or may enter the resistant state. No transition out of $R$ is modeled, so $R$ is absorbing under the forecast dynamics. All geographic visualizations use the Esri World Street Map basemap (Esri, DeLorme, NAVTEQ, USGS, Intermap, NRCAN, Esri Japan, METI, Esri China, and TomTom).

\begin{remark}[Reproducibility]
All the numerical results presented in this section can be reproduced using the code on GitHub at \href{https://github.com/Siming-Liang/Epidemic_ABM}{https://github.com/Siming-Liang/Epidemic\_ABM}.
\end{remark}

\subsection{Controlled validation of the assimilation framework}\label{sec:abm_uniform}

\begin{figure}[h!]
    \centering
\includegraphics[width=0.45\linewidth]{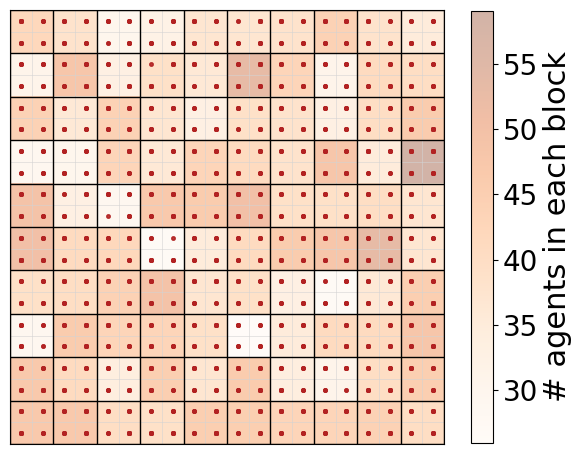}
    \caption{Uniform-grid ABM setup with $4000$ agents on a $20\times20$ grid. Red points show agent locations, while background shading gives the number of agents in each $2\times2$ observation block. The varying shading shows that agents are distributed nonuniformly across blocks.}
    \label{fig:abmgrid}
\end{figure}

To validate the proposed filtering framework in a controlled setting, we first consider an ABM evolving on a uniform spatial grid where the agents are distributed nonuniformly across the observation blocks, as illustrated in Figure~\ref{fig:abmgrid}. The domain consists of $20\times 20$ grid points (400 cells) populated by $J=4000$ agents initially distributed across the grid. Observations are not collected at the agent level; instead, we observe aggregated SIR counts over $2\times 2$ sub-grid blocks, resulting in $10\times 10=100$ observation blocks. Spatial interactions are local: disease transmission occurs only through contacts among agents in neighboring cells. To introduce heterogeneity and create a high-dimensional inference problem, we assign each observation block its own transmission and recovery parameters. Specifically, the transmission rate and recovery rate are sampled independently from a uniform distribution on $[0.1,0.3]$ for each observation block, resulting in $200$ unknown parameters in total. This configuration mimics the practical situation where epidemiological parameters vary across space and subpopulations \cite{bhattacharya_epihiper_2024}. We use an ensemble of size $L=100$. All ensemble members are initialized to match the same population counts in each observation block, while the initially infected agents are selected randomly (so the total number of initial infections is fixed but their locations vary across ensembles). Since the true parameters are assumed unknown, initial parameter guesses are sampled from $\mathcal{N}(0.4,0.1^2)$, intentionally biased away from the true range. Independent observation noise is added to the infected and resistant counts in each observation block, with a standard deviation equal to $10\%$ of the corresponding true count. In the experiment, we assimilate only $I$ and $R$. The susceptible count is then determined by local population conservation as $S=N-I-R$. To focus inference on transmission and recovery, the probability that an agent who clears infection enters $R$, termed the gain-resistance chance, is treated as known and sampled independently for each block from $\mathcal U[0.05,0.15]$.

\begin{figure}[h!]
    \centering
\includegraphics[width=0.85\linewidth]{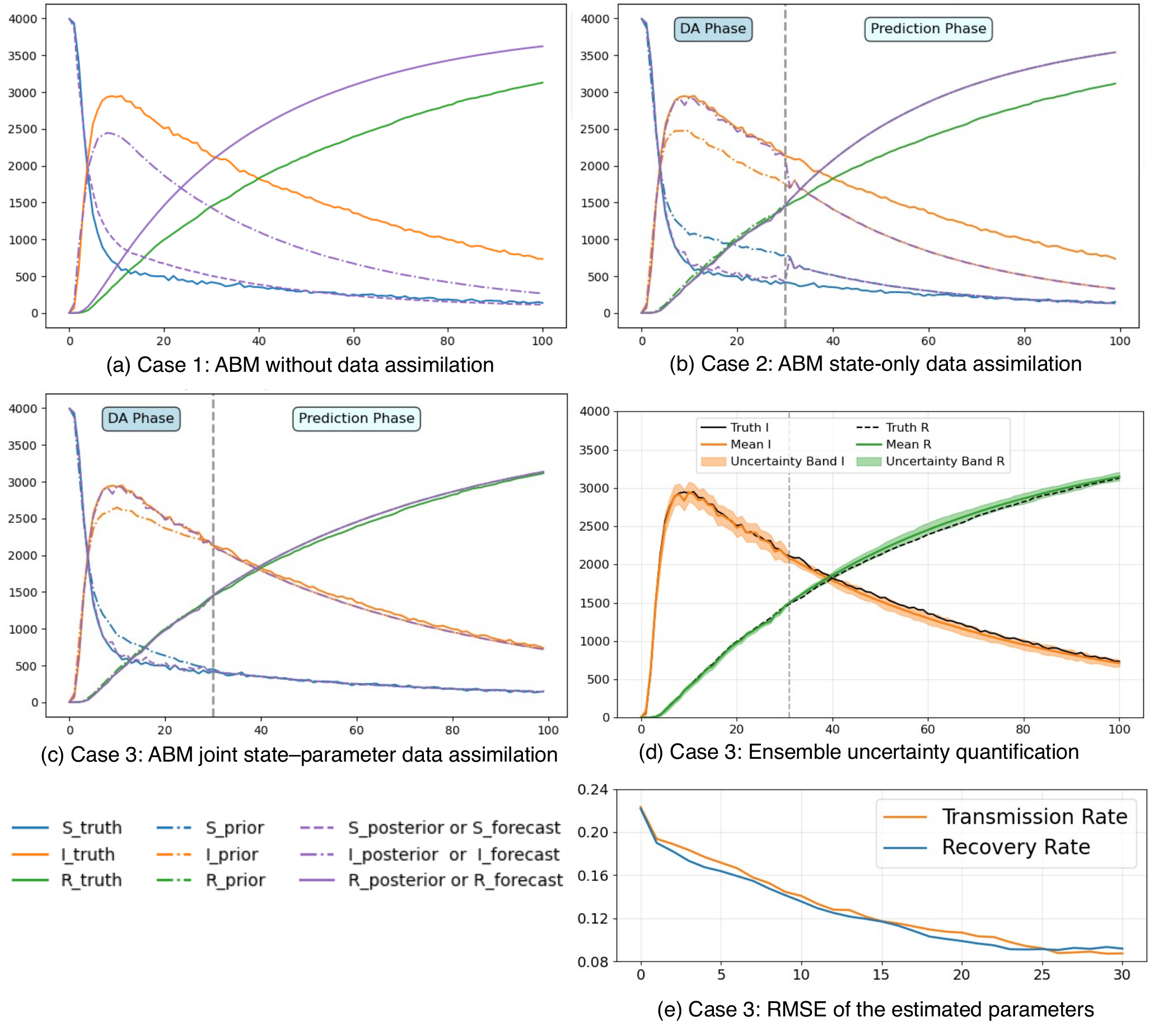}
    \caption{Uniform-grid SIR dynamics for (a) Case~1 without assimilation, (b) Case~2 with state-only assimilation, and (c) Case~3 with joint state--parameter assimilation. The vertical dashed line at $T=30$ marks the end of assimilation; panels (d) and (e) show Case~3 uncertainty and parameter RMSE, respectively. For either scalar parameter component, $\mathrm{RMSE}_{\theta}(t)=[K^{-1}\sum_{k=1}^{K}(\overline\theta_{t,k}-\theta_{k}^{\mathrm{true}})^2]^{1/2}$, where $\overline\theta_{t,k}=L^{-1}\sum_{l=1}^{L}\theta_{t,k,l}$, $K=100$, and $L=100$. Case~3 remains accurate during prediction, showing the importance of parameter learning.}
    \label{fig:case123anduq}
\end{figure}

For the experiment, we first simulate the ABM using the incorrect parameter ensemble without any data assimilation. Next, we apply data assimilation to the macrostate only during the first 30 time steps, while keeping parameters fixed at their incorrect initial values. Finally, we apply the proposed GenDA to jointly calibrate both the macrostate and the parameters during the first 30 time steps. Figure~\ref{fig:case123anduq} shows the aggregated SIR statistics over all observation blocks, i.e., the full-space dynamics.
\paragraph{Case 1: Uncorrected forecasts under misspecified transmission and recovery rates}
This baseline experiment establishes the effect of parameter misspecification in a stochastic ABM when no observational correction is applied. It isolates model error from data-assimilation effects and serves as a reference for evaluating whether subsequent improvements arise from state correction, parameter correction, or both. Figure~\ref{fig:case123anduq}(a) shows that the ensemble mean deviates substantially from the baseline (ground truth) trajectory. This demonstrates that, under parameter misspecification and stochastic dynamics, the ABM forecast alone is insufficient for accurate inference, motivating the need to incorporate observational information through data assimilation.

\paragraph{Case 2: Assimilating the epidemic state only}

This experiment tests whether correcting the epidemic state alone, while keeping incorrect parameters fixed, is sufficient for reliable prediction. Figure~\ref{fig:case123anduq}(b) shows that within the assimilation window ($t\le 30$), the ensemble estimate mean (posterior) corrects the model prediction (prior) mean and tracks the baseline closely, indicating that score-based nonlinear filtering effectively corrects the macrostate by combining model forecasts with observations. However, the underlying dynamics remain misspecified because the parameters are not updated. Once assimilation stops and the system enters a free prediction phase, the forecast deteriorates rapidly. This case illustrates the limitation of state-only assimilation and motivates joint state--parameter estimation.

\paragraph{Case 3: Jointly inferring epidemic state and parameters}
This experiment evaluates the central hypothesis of the paper: that parameters can be inferred indirectly by leveraging state corrections at the macro level, even when parameters are not directly observable. 
Figure~\ref{fig:case123anduq}(c) shows that the model prediction (prior) mean improves during assimilation and that the ensemble mean tracks the baseline closely. More importantly, the subsequent prediction phase remains accurate because joint state--parameter filtering modifies the model's internal dynamics rather than only correcting the current state. The overall decrease of the root mean square error (RMSE) in Figure~\ref{fig:case123anduq}(e) further quantifies effective online calibration of spatially heterogeneous parameters under the proposed GenDA. Together, these results show that joint state--parameter assimilation improves out-of-sample predictive performance by reducing model misspecification.

Joint state--parameter assimilation also reduces ensemble spread across the macrostates. During the assimilation window, uncertainty narrows as observations constrain both the epidemic state and the inferred parameters. Figure~\ref{fig:case123anduq}(d) visualizes this spread as a shaded ensemble envelope: the posterior tracks the baseline during assimilation, and the forecast remains relatively close to the baseline after assimilation ends.

\subsubsection{Block-level spatial accuracy and uncertainty}

\begin{figure}[!htb]
    \centering      \includegraphics[width=0.98\linewidth]{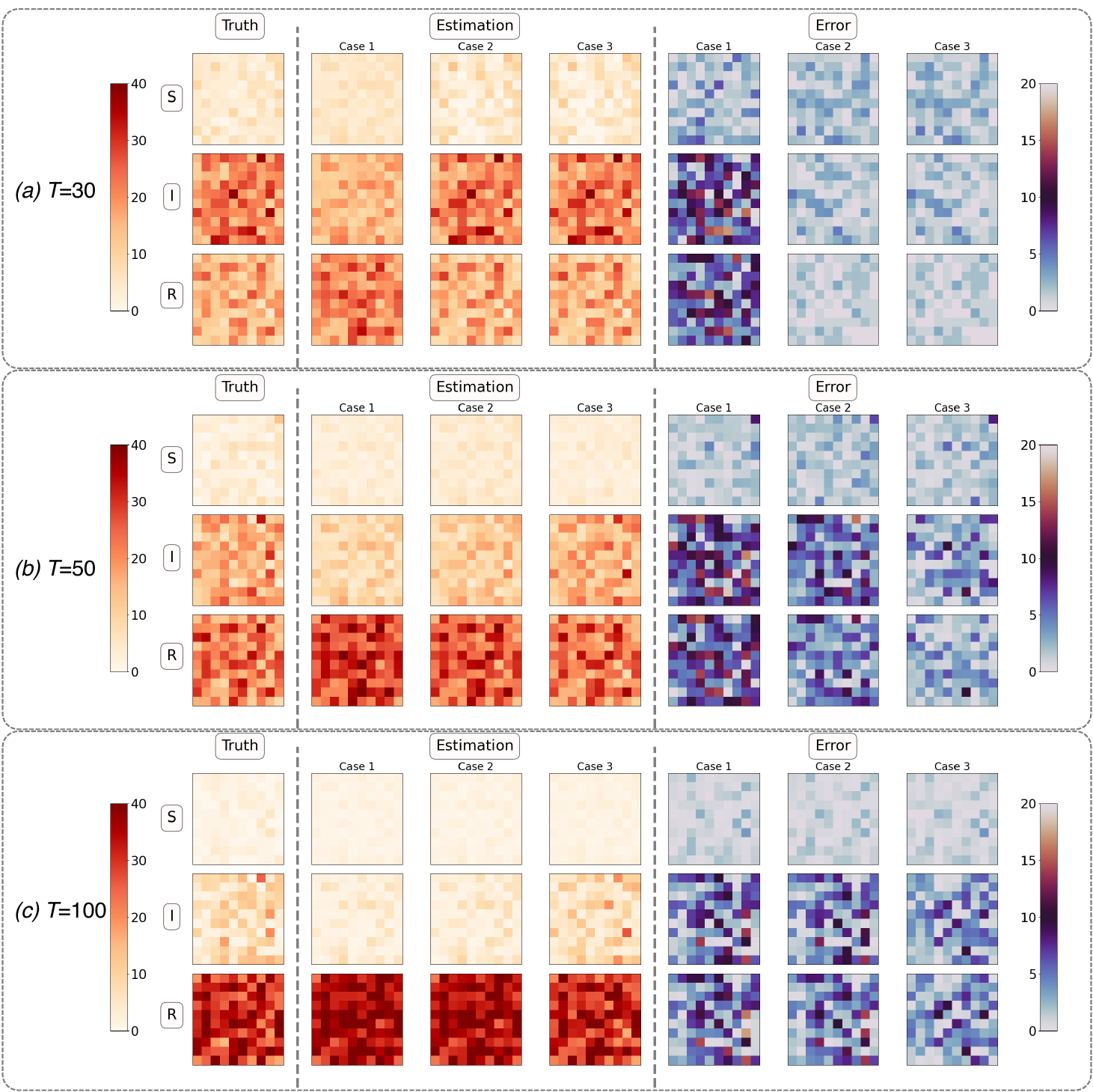}
    \caption{Uniform-grid SIR fields at $T=30$, $50$, and $100$. Rows show $S$, $I$, and $R$; columns show the synthetic reference, estimates from Case~1 (no assimilation), Case~2 (state-only assimilation), and Case~3 (joint state--parameter assimilation), followed by their absolute errors. Cases~2 and 3 are accurate at the end of assimilation ($T=30$), but Case~3 retains lower prediction errors after 20 and 70 prediction steps, highlighting the benefit of parameter correction for post-assimilation prediction. Spatial errors nevertheless grow at long forecast times because of stochastic agent dynamics.}

    \label{fig:grid_error_all}
\end{figure}

We next examine errors at the grid (macrostate) level for all three cases. Figure~\ref{fig:grid_error_all} compares spatial errors at $T=30$, $50$, and $100$, corresponding to the end of assimilation, a short prediction horizon, and a long prediction horizon, respectively. 

At $T=30$, Figure~\ref{fig:grid_error_all}(a) shows that, without data assimilation (Case 1), the ensemble does not reproduce the baseline spatial macrostate. In contrast, both Case 2 and Case 3 achieve low grid-level errors at the end of the assimilation window, consistent with the full-space statistics in Figures~\ref{fig:case123anduq}(b) and \ref{fig:case123anduq}(c). 
At $T=50$, after 20 steps of prediction, Figure~\ref{fig:grid_error_all}(b)  shows that the proposed GenDA (Case 3) still produces accurate macrostates with low spatial error. In contrast, Case 2 (state-only assimilation) deteriorates rapidly because the parameters remain incorrect, consistent with the full-space divergence observed in Figure~\ref{fig:case123anduq}(b).
At $T=100$, after 70 steps of free prediction, the inherent sensitivity of ABM dynamics becomes apparent. Figure~\ref{fig:grid_error_all}(c) shows that local spatial allocation errors grow even when full-space aggregated statistics remain accurate (Figure~\ref{fig:case123anduq}(c)). This discrepancy highlights a fundamental feature of stochastic ABMs: ensemble-based assimilation can recover accurate global statistics, but the {spatial allocation} of those totals remains sensitive to microscopic randomness such as mobility sampling, stochastic contacts, and discrete interaction events. Thus, full-space consistency does not guarantee long-range spatial accuracy.

\remark{Observation sparsity and scalability.}
The present experiments are conducted at moderate scale due to the computational limitations of \texttt{Mesa}, while the filtering methodology inherits the parallelism and scalability of the score-based filter, which has been demonstrated on large-scale HPC platforms in prior work \cite{rafid24scalable,sc24_1}.
In realistic surveillance settings, macro observations may be missing in parts of the domain. When observations are unavailable at certain locations or times, score filter based inpainting methods \cite{liang_ensf_inpainting_2025} and covariance-informed score filter variants \cite{zhang2025iensf, zhang2025exact} provide practical mechanisms to stabilize state estimation and 
support uncertainty characterization under sparse observation coverage.

\subsection{A geographically explicit epidemic case study}

\begin{figure}[!htb]
    \centering

    \begin{subfigure}{0.36\linewidth}
        \centering
        \includegraphics[width=\linewidth]{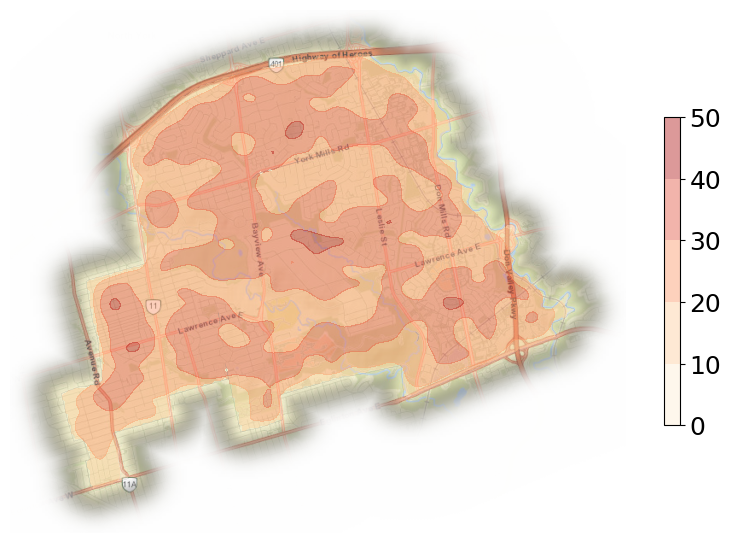}
        \caption{Agent density}
        \label{fig:agentdensity}
    \end{subfigure}
    \hfill
    \begin{subfigure}{0.31\linewidth}
        \centering
        \includegraphics[width=\linewidth]{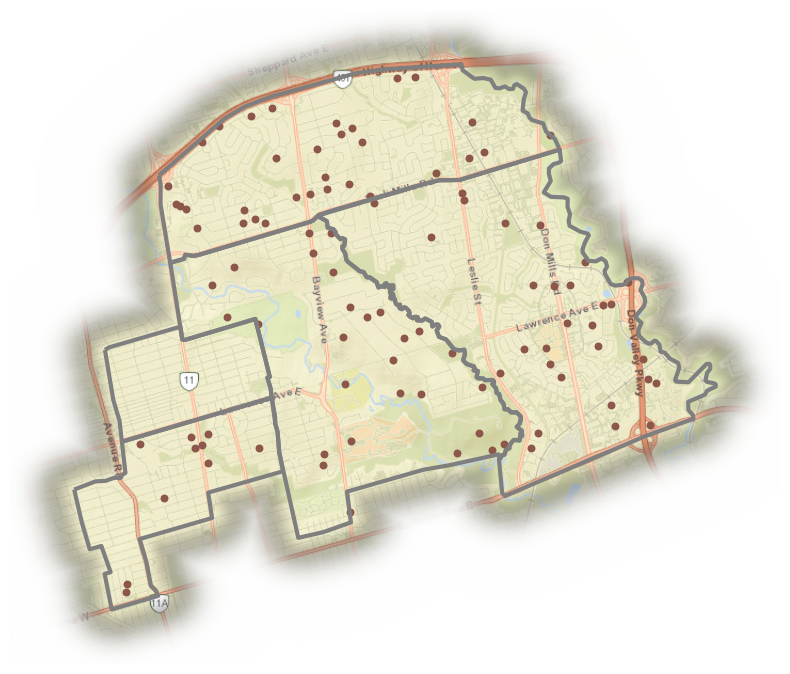}
        \caption{Observation locations}
        \label{fig:geo_grid}
    \end{subfigure}
    \hfill
    \begin{subfigure}{0.31\linewidth}
        \centering
        \includegraphics[width=\linewidth]{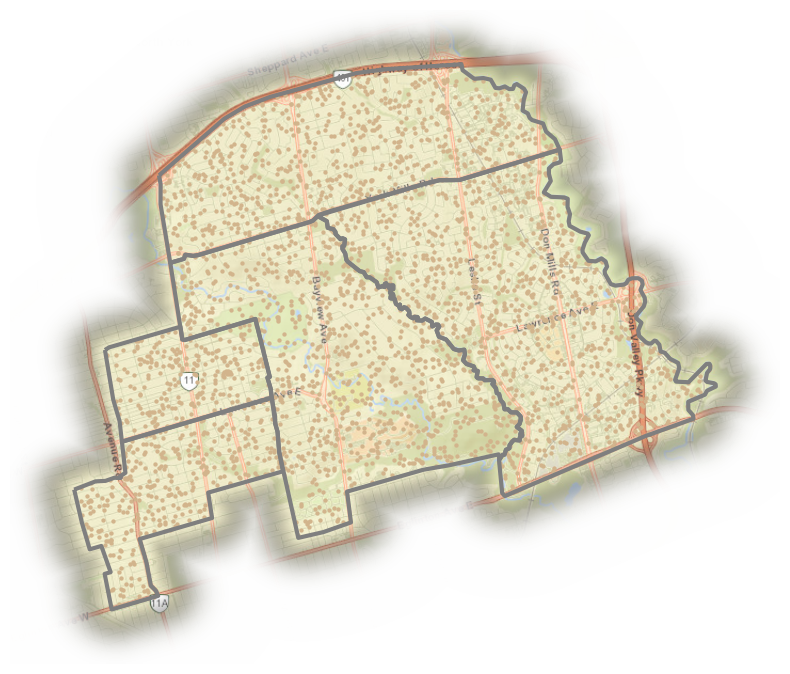}
        \caption{Hub locations}
        \label{fig:geo_hub}
    \end{subfigure}
    \caption{ABM-GEO configuration over five Toronto parameter subregions (see supplementary for details). (a) synthetic density of approximately $6000$ agents, (b) $100$ randomly sampled observation locations, and (c) approximately $3000$ simulated activity hubs. This setup tests GenDA under heterogeneous parameters, hub-driven mobility, and spatially uneven observation coverage.}
    \label{fig:combined_row}
\end{figure}

The uniform-grid experiments in Section~\ref{sec:abm_uniform} validate GenDA in a controlled setting, while the ABM-GEO experiments test the same framework in a geographically explicit, mobility-driven system with sparse observations. This setting is harder to calibrate because transmission depends on spatially varying population and hub density, as well as time-dependent agent activity and mobility. The hubs represent several activity types, including shopping centers, gyms, schools, workplaces, and residential locations. The goal is to assess whether GenDA retains state and parameter skill when the simulation includes geographically explicit geometric constraints.
We use a portion of the City of Toronto partitioned into five parameter subregions. The basemap and subregion geometries are geographic, all agents, hubs, epidemic trajectories, parameters, and surveillance observations are synthetic. Figure~\ref{fig:agentdensity} shows the simulated population density with approximately 6000 agents, and Figure~\ref{fig:geo_hub} shows roughly 3000 simulated hubs. To provide a controlled, nonoptimized sparse-observation layout, the \(100\) observation locations are sampled randomly and independently of the population-density field (Figure~\ref{fig:geo_grid}). The layout represents surveillance systems that provide limited and coarse-grained information. Observation sites near subregion boundaries may also reflect movement across administrative partitions, so the ABM-GEO setup combines spatially uneven observation with mobility-driven transmission. Individual agents move freely on the map according to an activity schedule updated every three data-assimilation cycles.

Following the experiment design in Section~\ref{sec:abm_uniform}, we again consider three cases: ABM-GEO forecast without data assimilation, state-only assimilation, and joint state--parameter assimilation. Transmission and recovery rates vary across the five subregions, with true values sampled from uniform distributions \([0.4,0.6]\) and \([0.1,0.3]\), respectively. We generate \(L=100\) ensembles with the same population counts at each observation location. The initial number of infected individuals is identical across ensembles, but infected agents are selected randomly. Since transmission and recovery rates are treated as unknown, initial parameter guesses are sampled from \(\mathcal{N}(0.3,0.1^2)\) for transmission and \(\mathcal{N}(0.4,0.1^2)\) for recovery. Observation noise has a standard deviation equal to \(10\%\) of the corresponding true macrostate count at each observation location.

\subsubsection{Regional epidemic burden over time}
The aggregate diagnostics for the ABM-GEO experiment are reported in Supplementary Figure~S1. These results are consistent with the uniform-grid experiment in Section~\ref{sec:abm_uniform}. State-only assimilation improves the observed epidemic trajectory but remains sensitive to parameter misspecification, whereas joint state--parameter assimilation reduces the parameter error and gives more stable aggregate epidemic estimates. We therefore focus the main-text figures on the spatial reconstruction behavior.

\subsubsection{Hotspot structure and regional error patterns}
\begin{figure}[!tb]
    \centering
    \includegraphics[width=0.98\linewidth]{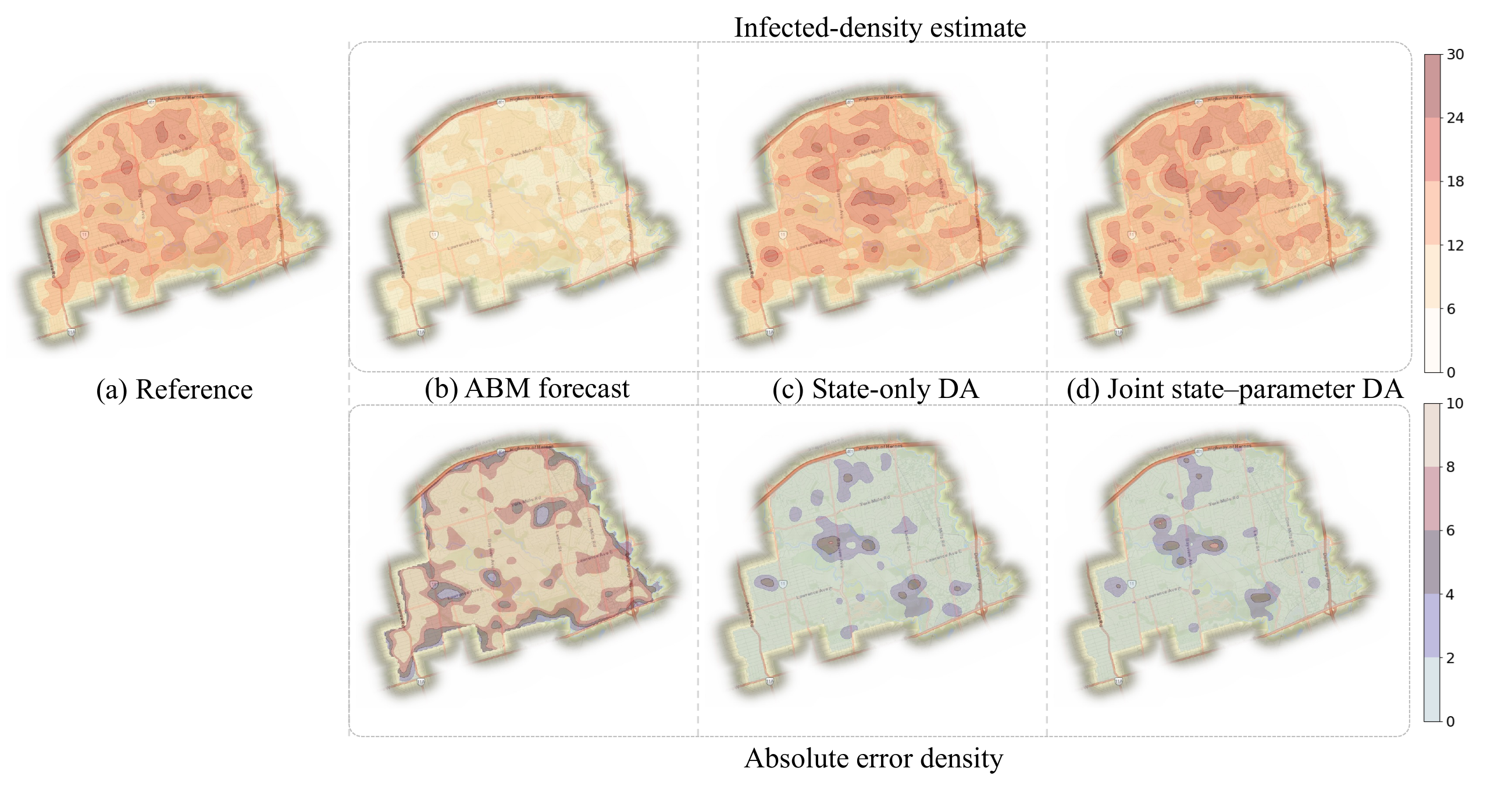}
    \caption{ABM-GEO infected-density fields at $T=30$, the end of assimilation. The top row shows the synthetic reference and estimates from Case~1 (no assimilation), Case~2 (state-only assimilation), and Case~3 (joint state--parameter assimilation); the bottom row shows their absolute errors. Cases~2 and 3 recover the dominant infection hotspots and substantially reduce the errors of Case~1, showing that aggregated observations can constrain the spatial epidemic state.}

\label{fig:all_cases_t30_with_colorbars}
\end{figure}

We now examine spatial behavior on the geographic map. Unlike the uniform-grid study (Figure~\ref{fig:grid_error_all}), we focus on the infected population because the magnitude and spatial distribution of infections are most relevant for epidemiological decision making. Figure~\ref{fig:all_cases_t30_with_colorbars} shows infected-density fields for all three cases, together with the reference solution and the corresponding error fields at \(T=30\). Because surveillance data constrain only aggregated macrostates, the reconstructed spatial fields should be interpreted as representative realizations rather than exact reconstructions. Except for the reference field, all estimated density fields are obtained through the macro--micro mapping procedure in Section~\ref{sec:Micro--macro}. The displayed spatial errors therefore reflect both macrostate discrepancies and additional variability introduced by stochastic reassignment, which is inherent to agent-based representations.

At \(T=30\), Case~1 (ABM-GEO forecast without data assimilation) fails to reproduce the spatial structure of the infected-density field and exhibits large, spatially coherent errors across the domain (Figure~\ref{fig:all_cases_t30_with_colorbars}(b)). In contrast, both Case~2 (state-only assimilation) and Case~3 (joint state--parameter assimilation) produce infected-density fields that are much closer to the baseline solution, with substantially reduced error magnitudes (Figure~\ref{fig:all_cases_t30_with_colorbars}(c)--(d)). Despite sparse observations and stochastic agent dynamics, the score-based filter constrains the ABM-GEO macrostate, suppresses large-scale spatial biases, and recovers the dominant infection patterns.

\begin{figure}[!htb]
    \centering
    \includegraphics[width=0.99\linewidth]{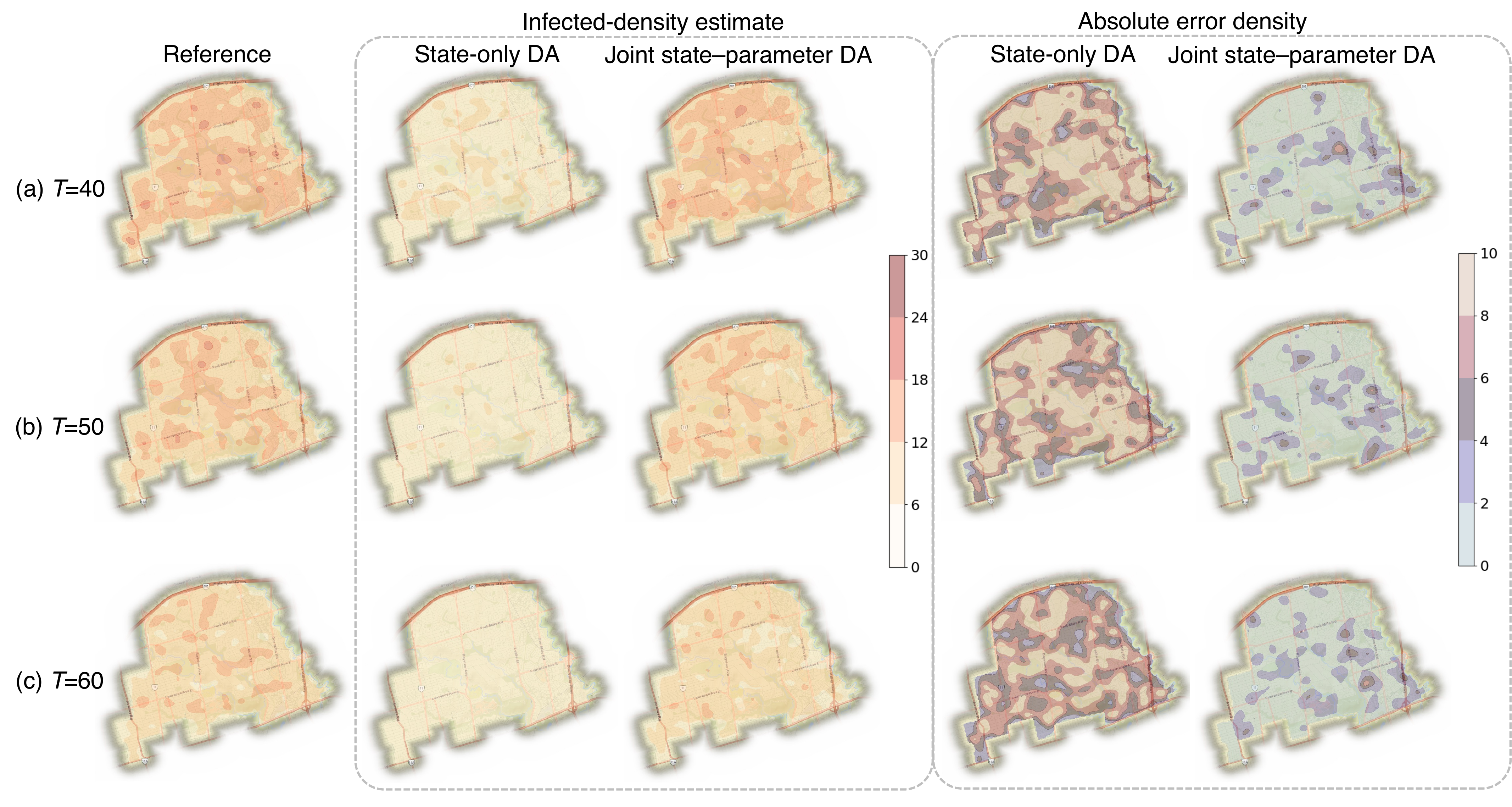}
    \caption{ABM-GEO infected-density predictions at $T=40$, $50$, and $60$, after assimilation ends at $T=30$. Columns show the synthetic reference, Case~2 state-only prediction, Case~3 joint state--parameter prediction, and their absolute errors. Case~2 progressively loses the hotspot structure, whereas Case~3 maintains smaller errors and the dominant regional patterns, demonstrating the importance of parameter learning for post-assimilation prediction.}
\label{fig:geo_error_case2_case3}
\end{figure}

Figure~\ref{fig:geo_error_case2_case3} compares the ABM-GEO infected-density estimates and absolute error densities for Case~2 and Case~3 at $T=40$, $50$, and $60$. The comparison highlights the effect of parameter learning during the prediction phase. In Case~2, state-only assimilation is vulnerable to parameter misspecification: the predicted infection density becomes systematically weaker than the reference field, dominant high-density regions fade, and spatially coherent error structures persist across the forecast window. In contrast, Case~3 substantially reduces these large-scale errors. The joint state--parameter update preserves the main regional infection patterns, including the broad spatial extent of high-risk regions and the elongated hotspot structures visible in the reference solution.
The remaining discrepancies in Case~3 are expected in a stochastic ABM-GEO system with sparse and spatially uneven observations. The ensemble-mean density represents an average over plausible epidemic trajectories, so local peaks may be smoothed and fine-scale boundaries may differ from the reference realization. The relevant diagnostic is therefore whether the method retains the dominant hotspot structure and regional risk gradients. From this perspective, Case~3 provides meaningful spatial skill, whereas Case~2 exhibits parameter-driven spatial drift after the assimilation window.

%% file: 5Conclusion.tex
\section{Discussion and Conclusion} \label{sec:conclusion}

This paper develops the GenDA framework for a setting with sparse, noisy, and aggregated data and a stochastic agent-based simulator as the forward model. The framework operates at the level of epidemic macrostates while preserving consistency with agent-based dynamics. Its goal is not to reconstruct an unobservable ``true'' agent configuration, but to recover macro-scale epidemic structure that can support interpretation and scenario analysis.
The numerical results illustrate this point in two complementary synthetic settings. In the uniform-grid experiment, state-only assimilation improves short-term tracking but loses predictive skill once observations are no longer incorporated, whereas joint state--parameter filtering stabilizes both the inferred epidemic state and the subsequent forecast. In the ABM-GEO experiment, the same pattern persists under uneven population density, spatially uneven observation support, stochastic mobility, and hub-mediated contact structure. GenDA preserves support-level consistency and recovers burden trends and dominant hotspot structures after the assimilation window.

For epidemic surveillance, these findings clarify what can reasonably be expected from a calibrated epidemic ABM. Such a model should not be judged by its ability to reproduce one realized agent trajectory in full detail. More appropriate targets are support-level burdens, spatial risk gradients, and forecast distributions. Biologically, the inferred transmission and recovery fields should therefore be interpreted as effective parameters that summarize heterogeneous contact and clearance processes at the chosen sharing scale, rather than as uniquely identified individual-level mechanisms.

Several limitations define the scope of the evidence. Both experiments are synthetic, and the Toronto case uses real geometry but not real population, facility, or surveillance data; empirical validation remains necessary before operational use. The experiments are moderate in scale and do not establish end-to-end high-performance computing scalability. Parameter identifiability may weaken when distinct parameter fields produce similar aggregate trajectories, particularly under sparse supports, model-form error, or misspecified observation noise. The exact local reassignment used here also assumes disjoint supports; overlapping supports require a joint constrained allocation. Finally, the displayed ensemble envelopes summarize spread for single synthetic reference trajectories and do not constitute a repeated-sampling coverage study.

A primary next step is integration with the ENABLE high-performance agent-based framework \cite{spannaus_enable_2025}, which is designed for state- and national-scale population simulations using real-world data. That integration would make it possible to test the same inferential ideas in substantially larger epidemic systems. Further work should also evaluate sensitivity to support geometry and reporting error, address overlapping constraints, and validate parameter and state estimates against empirical surveillance.

Overall, the results show that partially observed epidemic ABMs can be calibrated coherently at the aggregate level. By combining non-Gaussian macrostate assimilation, direct parameter updates, and macro--micro consistency, GenDA extracts useful epidemic summaries from systems in which exact fine-scale predictability is fundamentally limited. The framework provides a testable foundation for future work with larger simulators and empirical surveillance data.

%% file: 6Appendix.tex
\appendix
\section{A Reader's Guide for Life-Science Audiences}

\subsection{Interpreting the filtering cycle}

The agent-based model used in this paper evolves at the level of individual agents, but the available data are regional surveillance summaries. Therefore, the method is not trying to reconstruct the exact hidden history of every agent. Many different agent-level configurations can produce the same regional counts, especially in a stochastic epidemic model.
Instead, the inferred quantities are \emph{macrostates}: sub-region level epidemic counts and a heterogeneous parameter field that summarizes transmission and recovery conditions. These are also the quantities most relevant for public-health interpretation, such as regional burden, hotspot persistence, and forecast uncertainty.
Each filtering cycle can be read as a repeated forecast--correction step. First, the ABM is advanced to produce a forecast. Second, surveillance data are used to correct the sub-regionlevel macrostate and parameter field. Third, the corrected macrostate is mapped back to a plausible agent-level configuration so that the ABM can continue running. This last step should be understood as a consistency step, not as a unique reconstruction of individual-level trajectories.

\subsection{Further reading}

For readers interested in the biological and modeling background of epidemic ABMs, useful starting points include the general ABM texts of Railsback and Grimm \cite{railsback_grimm_2019}, recent reviews of epidemic ABMs \cite{zhang2025agent}, and large-scale epidemic modeling platforms such as EpiPredict \cite{suer_epipredict_2024}, UVA-EpiHiper \cite{bhattacharya_epihiper_2024}, and ENABLE \cite{spannaus_enable_2025}. For readers interested in uncertainty quantification and calibration for stochastic ABMs, relevant entry points include \cite{gugole_uq_covid_abm_2021,kimpton_uq_abm_2024,swallow_emulation_abm_2022}. For readers interested in data assimilation for ABMs and related biomedical applications, helpful references include \cite{cocucci_enkf_epi_abm_2022,knapp2025personalizing,ward_enkf_abm_2016}. Finally, the score-based filtering perspective used in this paper builds on the recent score-based filter literature \cite{bao2025nonlinear,bao2024ensemble,bao2024score,bao2023unified}.